\documentclass{amsart}
\usepackage{latexsym}
\usepackage[all]{xy}

\usepackage{amssymb} % \mathbb
\usepackage{amsmath} % \mathbb
\usepackage{amsthm}

\def\Z{{\mathbb{Z}}}% \Z == \mathbb{Z}
\def\K{{\mathbb{K}}}% \K == \mathbb{K}
\def\A{{\mathcal{A}}}% \A == \mathcal{A}
\numberwithin{equation}{section}

\newcommand{\owari}{\hfill$\square$}

\newtheorem{theorem}{Theorem}[section]
\newtheorem{prop}[theorem]{Proposition}
\newtheorem{cor}[theorem]{Corollary}
\newtheorem{lemma}[theorem]{Lemma}

\newtheorem{conjecture}[theorem]{Conjecture}

\newcommand{\xgraphAvertex}[1][****]{
\xgraphAVertex #1
}
\newcommand{\xgraphAVertex}[4]{
\if#3o\put(0,0){\circle{4}}\fi%
\if#3*\put(0,0){\circle*{4}}\fi%
\if#3.\put(0,0){\circle*{2.4}}\fi%
\if#4o\put(30,0){\circle{4}}\fi%
\if#4*\put(30,0){\circle*{4}}\fi%
\if#4.\put(30,0){\circle*{2.4}}\fi%
\if#2o\put(0,30){\circle{4}}\fi%
\if#2*\put(0,30){\circle*{4}}\fi%
\if#2.\put(0,30){\circle*{2.4}}\fi%
\if#1o\put(30,30){\circle{4}}\fi%
\if#1*\put(30,30){\circle*{4}}\fi%
\if#1.\put(30,30){\circle*{2.4}}\fi%
}
\newcommand{\xgraphA}[6]{
\if#1+\put(30,30){\line(-1,0){30}}\fi% ~
\if#1.\qbezier[7](30,30)(15,30)(0,30)\fi% ~
\if#1-\put(30,31){\line(-1,0){30}}\put(30,29){\line(-1,0){30}}\fi% ~%%
\if#2+\put(0,0){\line(1,1){30}}\fi % /
\if#2.\qbezier[10](0,0)(15,15)(30,30)\fi % /
\if#2-\put(-0.7,0.7){\line(1,1){30}}\put(0.7,-0.7){\line(1,1){30}}\fi % /%%
\if#3+\put(30,30){\line(0,-1){30}}\fi% )
\if#3.\qbezier[7](30,0)(30,15)(30,30)\fi % /
\if#3-\put(29,30){\line(0,-1){30}}\put(31,30){\line(0,-1){30}}\fi% )%%
\if#4+\put(0,0){\line(0,1){30}}\fi % (
\if#4.\qbezier[7](0,0)(0,15)(0,30)\fi % /
\if#4-\put(-1,0){\line(0,1){30}}\put(1,0){\line(0,1){30}}\fi % (
\if#5+\put(0,30){\line(1,-1){30}}\fi % \
\if#5.\qbezier[10](30,0)(15,15)(0,30)\fi % /
\if#5-\put(-0.7,29.3){\line(1,-1){30}}\put(0.7,30.7){\line(1,-1){30}}\fi % \%%
\if#6+\put(0,0){\line(1,0){30}}\fi %_
\if#6.\qbezier[7](0,0)(15,0)(30,0)\fi % /
\if#6-\put(0,1){\line(1,0){30}}\put(0,-1){\line(1,0){30}}\fi %_
\xgraphAvertex}
\begin{document}

\title[On the conjecture of Athanasiadis]
{On a conjecture of Athanasiadis
related to freeness of a family of hyparplane arrangements
% : Appendix to Signed-eliminable graphs and free multiplicities on the braid arrangement
}
\author{Takuro Abe}
\address{Department of Mechanical Engineering and Science, Kyoto University, Yoshida Hon-Machi, 
Sakyo-Ku, Kyoto 606--8501, Japan}
\email{abe.takuro.4c@kyoto-u.ac.jp}
%\author{Koji Nuida}
%\address{Research Center for Information Security (RCIS), National Institute of
%Advanced Industrial Science and Technology (AIST),
%Akihabara-Daibiru Room 1003, 1-18-13 Sotokanda, Chiyoda-ku, Tokyo
%101--0021, Japan}
%\email{k.nuida@aist.go.jp}
%\author{Yasuhide Numata}
%\address{Faculty of Integrated Media, Wakkanai Hokusei Gakuen University, 
%Wakabadai 1-2290-28, Wakkanai, Hokkaido 097--0013, Japan}
%\email{numata-y@wakhok.ac.jp}



%\pagestyle{plain}

%%%%%%    TEXT START    %%%%%%

\maketitle

\begin{abstract}
We prove a characterization
of freeness, conjectured by Athanasiadis, for the family of
hyperplane arrangements which lie between the Coxeter and the
Catalan arrangement of type $A_\ell$. One direction was already
proved in \cite{ANN}. Here we prove the other direction
%We prove 
%a characterization of freeness, 
%conjectured by Athanasiadis, for deformations of the braid
%arrangement which lie between the Coxeter and Catalan arrangement of 
%type $A_\ell$ 
%in terms of some digraphs. 
%the conjecture of Athanasiadis, which states a characterization of 
%the freeness of a deformation of the braid arrangement between the Coxeter and Catalan arrangement of 
%type $A_\ell$ 
%in terms of some digraphs. 
%One direction of it was already proved in \cite{ANN}, and we prove the 
%converse direction.
\end{abstract}
\setcounter{section}{-1}

\section{Conjecture of Athanasiadis}
The problem to characterize freeness has been an
important and interesting one within the theory of hyperplane
arrangements. 
%In the study of arrangements of hyperplanes, to characterize their freeness
%has been an important problem and there are a lot of interesting 
%researches on it. 
For example, the freeness of graphic arrangements
%, i.e., 
%a subarrangement of the braid arrangement, 
is characterized by Stanley 
in terms of chordal graphs (see 
\cite{S} and \cite{ER} for details). Also, 
free arrangements between the Coxeter arrangements of the types $A_{\ell-1}$ and $B_\ell$ 
are characterized by Edelman and Reiner in \cite{ER}, 
and those between the Coxeter and Shi arrangements of type $A_\ell$ by Athanasiadis in \cite{Ath1}. 
A conjecture 
of Athanasiadis, which was introduced in \cite{Ath2}, 
is a generalization of his freeness characterization in \cite{Ath1}.
%also related to a characterization of the freeness of some arrangements. 
To state the conjecture, let us recall several definitions and 
results. 
In this article we use \cite{OT} as a general reference, and the same notation as in \cite{ANN}. 

Let $\K$ be an arbitrary field of characteristic zero 
and consider an affine arrangement in $V^{\ell+1}=\K^{\ell+1}$ defined by 
\begin{eqnarray}\label{defath}
x_i-x_j=-k-\epsilon(i,j),-k,-(k-1),\ldots,k,k+\epsilon(j,i)  \\
 (1 \le i<j \le \ell+1), \nonumber
\end{eqnarray}
where $k \in \Z_{\ge 0}$ and $\epsilon(i,j)=0$ or $1$. Note that 
we distinguish $(i,j)$ and $(j,i)$ as explained later. 
Such arrangements are examples of 
deformations of the braid arrangement. 
%The main focus of these authors is on the characteristic polynomial of these arrangements. 
%Now let us go back to the deformation (\ref{defath}).
To consider the arrangement above, Athanasiadis introduced 
the directed graph $G$ consisting of the vertex set 
$V_G=\{1,2,\ldots,\ell+1\}$ and the set of directed edges 
$E_G \subset \{(i,j)|1 \le i,\ j \le \ell+1,\ i \neq j\}$. Here the edge 
$(i,j)$ is the arrow from $i$ to $j$. If we define 
\[
\epsilon(i,j):=
\left\{
\begin{array}{rl}
1 & \mbox{if}\ (i,j) \in E_G ,\\
0 & \mbox{if}\ (i,j) \not \in E_G, 
\end{array}
\right.
\]
then there is a one to one correspondence between 
the deformations of the form (\ref{defath}) and the digraph above.
%every affine arrangement 
%above can be expressed by these directed graphs. For such a graph $G$ let $\A_G$ denote 
%the corresponding arrangement of the form (\ref{defath}). 
In \cite{Ath0}, Athanasiadis gave 
a product formula of the characteristic polynomial of $\A_G$ when $G$ satisfies 
the following two conditions: 

\begin{itemize}
\item[(A1)] 
For every triple $i,j,k$ with $i,j<k$, it holds that, if $(i,j) \in E_G$, then 
$(i,k) \in E_G$ or $(k,j) \in E_G$.
\item[(A2)] For every triple $i,j,k$ with $i,j<k$, it holds that, 
if $(i,k) \in E_G$ and $(k,j) \in E_G$ then $(i,j) \in E_G$. 
\end{itemize} 

%Let us call these two conditions \textbf{Athanasiadis' conditions}.
%, and \textbf{A-conditions} for 
%short. 
%Also, we say a digraph 
%$G=(V_G,E_G)$ is an \textbf{Athanasiadis' digraph} if it satisfies 
%Athanasiadis' conditions. 
Note that conditions (A1) and (A2) depend on the ordering of the vertices 
$V_G$. Hence we say a digraph $G$ satisfies (A1) and (A2) if the graph satisfies them  
after re-ordering the set of vertices
$V_G$.
%=\{1,2,\ldots,\ell+1\}$. 
Based on these results and definitions, 
Athanasiadis also gave the following conjecture. 

\begin{conjecture}[\cite{Ath2}, Conjecture 6.6]\label{conjath}
Let $k=0$ in the deformation 
(\ref{defath}). Then the coning 
$c\A_G$ of $\A_G$ is free if and only if $G$ satisfies (A1) and (A2). 
\end{conjecture}

In \cite{ANN}, it was proved that 
conditions (A1) and (A2) are sufficient for Conjecture \ref{conjath} 
in more general setting, i.e., the following holds.

\begin{theorem}[\cite{ANN}, Theorem 5.3]
In the deformation (\ref{defath}), 
$c\A_G$ is free if $G$ satisfies (A1) and (A2). 
%In particular, the ``if'' part of Conjecture \ref{conjath} is true. 
\label{suff}
\end{theorem}

In this article we prove the converse of Theorem \ref{suff} as follows:

\begin{theorem}
In the deformation (\ref{defath}), 
if $c\A_G$ is free, then $G$ satisfies (A1) and (A2). 
%In particular, the ``only if'' part of Conjecture \ref{conjath} is true.
%, hence 
%Conjecture \ref{conjath} is true. 
\label{nece}
\end{theorem}

As a corollary we can show the conjecture of Athanasiadis.

\begin{cor}
The Conjecture \ref{conjath} of Athanasiadis is true.
\end{cor}

In the next section we prove Theorem \ref{nece}. The idea of the proof is 
to lift up the signed eliminable ordering of the signed graph to 
that on the digraph which makes $G$ satisfy (A1) and (A2). For the proof, we give 
a characterization of digraphs satisfying (A1) and (A2) independent of the numbering of 
vertices. 
\medskip

\noindent
\textbf{Acknowledgements}. 
The author is grateful to Christos Athanasiadis for helpful comments to this article. 
The author is partially supported 
by JSPS Grants-in-Aid for Young Scientists
(B) No. 21740014.

\section{Proof of Theorem \ref{nece}}

To prove Theorem \ref{nece} we need two results. First is the following:

\begin{lemma}
If a digraph $G$ contains one of the following induced subgraphs $H \subset G$ with 
$V_H=\{i,j,k\}$, then 
$c\A_G$ is not free:
\begin{itemize}
\item[(1)] 
%$V_H=\{i,j,k\}$ and 
$E_H=\{(i,j),(j,k)\}$.
\item[(2)]
%$V_H=\{i,j,k\}$ and 
$E_H=\{(i,j),(j,k),(k,i)\}$.
\item[(3)]
%$V_H=\{i,j,k\}$ and 
$E_H=\{(i,j),(j,k),(k,i),(i,k)\}$.
\end{itemize}
\label{lemma1}
\end{lemma}

\noindent
\textbf{Proof}. Let $H_\infty$ denote the infinite hyperplane of $c\A_G$ 
added to $\A_G$. 
If $c\A_G$ is free, then every localization of it is also free. 
Also, if $H$ is an induced subgraph of $G$ corresponding to the set of 
vertices $V_H=\{i,j,k\} \subset V_G$, then $c\A_H \times \emptyset_{\ell-2}=(c\A_G)_X$, where 
$X=\{x_i=x_j=x_k\} \cap H_\infty$ and $\emptyset_{\ell-2}$ is an $(\ell-2)$-dimensional 
empty arrangement. Since the freeness of $c\A_H \times \emptyset_{\ell-2}$ is equivalent to 
that of $c\A_H$, 
we may assume that $\ell=2,i=1,j=2,k=3$ and check the non-freeness of 
the deformation of $A_2$-type arrangements. For that purpose, 
we check the characteristic polynomial of (1), (2) and (3). We can obtain them by 
\cite{A1}, or direct computations. For (1), 
$\chi(\A_G,t)=t(t^2-(6k+5)t+(9k^2+15k+7))$. 
For (2), 
$\chi(\A_G,t)=t(t^2-(6k+6)t+(9k^2+18k+11))$. 
For (3), 
$\chi(\A_G,t)=t(t^2-(6k+7)t+(9k^2+21k+13))$. Since they are irreducible over $\Z$, 
Terao's factorization theorem in \cite{T} completes the proof. \owari
\medskip

Second result is a characterization of digraphs satisfying (A1) and 
(A2) without using the numbering of vertices. To state it, let us introduce a 
correspondence map $S$ from a digraph $G=(V_G,E_G)$ to a 
signed graph $S(G):=\overline{G}=(V_{\overline{G}},E_+ \cup E_-)$ with 
$E_\mu \subset \{\{i,j\} \mid 1 \le i,j \le \ell+1,\ i \neq j,\ 
\{i,j\}=\{j,i\}\}\ (\mu\in \{+,-\})$ and 
$E_+\cap E_-=\emptyset$. 
First, the sets of vertices are the same; 
$V_G=V_{\overline{G}}$. The relation between edges are as follows:
\begin{itemize}
\item
If $(i,j) \in E_G$ and $(j,i) \in E_G$, then 
$\{i,j\} \in E_+$.
\item 
If exactly one of $(i, j)$ and $(j, i)$ 
belongs to $E_G$, then 
$\{i,j\} \not \in E_+ \cup E_-$. 
\item
If $(i,j) \not \in E_G$ and $(j,i) \not \in E_G$, then 
$\{i,j\} \in E_-$. 
\end{itemize}

Also, recall the definition of a \textbf{signed eliminable} graph 
introduced in \cite{ANN}. We say a signed graph $(V_{\overline{G}},E_+\cup E_-)$ is signed 
eliminable if there is a numbering $\{1,2,\ldots,\ell+1\}$ of $V_{\overline{G}}$ such that 
for any $i,j,k \in V_{\overline{G}}$ with $i,j <k$, the following two conditions are satisfied:

\begin{itemize}
\item[(SE1)]
If $\{k,i\} \in E_\mu$ and $\{i,j\} \in E_\nu$ for $\{\mu,\nu\}=\{+,-\}$, then 
$\{k,j\} \in E_\nu$.
\item[(SE2)]
If $\{k,i\} \in E_\mu$ and $\{k,j\} \in E_\mu$ for $\mu \in \{+,-\}$, then 
$\{i,j\} \in E_\mu$. 
\end{itemize}

For details, see \cite{ANN} and \cite{N}. Note that a characterization of 
a signed eliminable graph without using a numbering of vertices is given in \cite{N}.

\begin{prop}
Let $G$ be a digraph and $S(G)=\overline{G}$ a signed graph defined above. 
Then $G$ satisfies (A1) and (A2) if and only if $\overline{G}$ 
is signed eliminable and $G$ does not contain any induced 
subgraph of three-vertices (1),(2) and (3) in Lemma \ref{lemma1}.
\label{char}
\end{prop}

\noindent
\textbf{Proof}.
First assume that $G$ satisfies (A1) and (A2). Then 
it was proved in \cite{ANN} that the corresponding signed graph $S(G)=\overline{G}$ is signed eliminable. 
Also, by definition, any induced subgraph of three vertices of $G$ 
satisfies (A1) and (A2). Since any numbering on vertices of digraphs (1), (2) and (3) in Lemma \ref{lemma1} 
cannot make them satisfying (A1) and (A2), $G$ cannot contain any of these subgraphs. 

Next assume that $\overline{G}$ is signed eliminable and $G$ does not contain 
any of (1),(2) and (3) in Lemma \ref{lemma1}. 
We may assume that $1,2,\ldots,\ell+1$ is a signed elimination ordering. We claim that 
the same ordering on $V_{\overline{G}}=V_G$ makes $G$ into a digraph satisfying 
(A1) and (A2), which is what 
we want to prove. 

To check this claim, it suffices to check that, for an 
induced subgraph $\overline{H}=\{i,j,k\} \subset V_{\overline{G}}$ 
of $\overline{G}$ with 
$i,j <k$, the induced subgraph $H:=\{i,j,k\} \subset V_G$ of $G$ satisfies two conditions (A1) and (A2). 
%Since the set of edges of ${\overline{H}}$ is the restriction of $E_{^{\overline{G}}}$ 
%onto ${\overline{H}}$, 
Let us denote $\overline{H}=\{\{i,j,k\},E_{+}|_{{\overline{H}}} \cup E_{-}|_{{\overline{H}}})$. 
This lifting correspondence from $\overline{H}$ to $H$ is as follows:
\begin{itemize}
\item
If $\{i,j\} \in E_+|_{{\overline{H}}}$, then $(i,j) \in E_G$ and $(j,i) \in E_G$.
\item
If $\{i,j\} \not \in E_+|_{{\overline{H}}} \cup E_-|_{{\overline{H}}}$, then either $(i,j) \in E_G$ or $(j,i) \in E_G$.
\item
If $\{i,j\} \in E_-|_{{\overline{H}}}$, then $(i,j) \not \in E_G$ and $(j,i) \not \in E_G$. 
\end{itemize}
So there are only finite possibility of the induced subgraph 
$H$ of $G$ as a lifting of 
$\overline{H}$. 
Let us check the lifting $H$ of $\overline{H}$ in every case. 
Put $E|_{{\overline{H}}}:=E_+|_{{\overline{H}}} \cup E_-|_{{\overline{H}}}$.
\medskip

\textbf{Case 1}. $E|_{{\overline{H}}}=\emptyset$. 
In this case, using symmetry, we have four possibilities of $H$:

\begin{enumerate}
\item
$\{(i,j),\ (i,k),\ (k,j) \}=E_H$.
\item
$\{(i,j),\ (i,k),\ (j,k)\} = E_H$.
\item
$\{(i,j),\ (k,i),\ (k,j) \}= E_H$.
\item
$\{(i,j),\ (k,i),\ (j,k) \}= E_H$.
\end{enumerate}

In these four cases, only (4) does not satisfy (A1) and (A2), which 
cannot occur by assumption.
\medskip

\textbf{Case 2}. $|E_+|_{{\overline{H}}}|=1$ and $E_-|_{{\overline{H}}}=\emptyset$. 
In this case, using symmetry, we have seven possibilities of $H$:

\begin{enumerate}
\item
$\{(i,j),\ (j,i),\ (i,k),\ (k,j) \}= E_H$.
\item
$\{(i,j),\ (j,i), \ (i,k),\ (j,k) \}= E_H$.
\item
$\{(i,j),\ (j,i),\  (k,i),\ (k,j) \}=E_H$.
\item
$\{(i,j),\ (k,i),\ (i,k),\ (k,j) \}= E_H$.
\item
$\{(j,i),\ (k,i),\ (i,k),\ (k,j) \}= E_H$.
\item
$\{(i,j),\ (k,i),\ (i,k),\ (j,k) \}= E_H$.
\item
$\{(j,i),\ (k,i),\ (i,k),\ (j,k) \}= E_H$.
\end{enumerate}

In these seven cases,  (1), (5) and (6) do not satisfy (A1) and (A2), which 
cannot occur by assumption.
\medskip

\textbf{Case 3}. $E_+|_{{\overline{H}}}=\emptyset$ and $|E_-|_{{\overline{H}}}|=1$. 
In this case, using symmetry, we have seven possibilities of $H$:

\begin{enumerate}
\item
$\{(i,k),\ (k,j) \}= E_H$.
\item
$\{(i,k),\ (j,k) \}= E_H$.
\item
$\{(k,i),\ (k,j) \}= E_H$.
\item
$\{(i,j),\ (j,k) \}= E_H$.
\item
$\{(i,j),\ (k,j) \}= E_H$.
\item
$\{(j,i),\ (j,k) \}= E_H$.
\item
$\{(j,i),\ (k,j) \}= E_H$.
\end{enumerate}

In these seven cases,  (1), (4) and (7) do not satisfy (A1) and (A2), which 
cannot occur by assumption.
\medskip

\textbf{Case 4}. $|E_+|_{{\overline{H}}}|=2$ and $E_-|_{{\overline{H}}}=\emptyset$. 
In this case, 
recall that $\overline{G}$ is, and so $\overline{H}$ is signed eliminable. Hence 
we can restrict the possibility of $H$ and $\overline{H}$. For example, 
$$
\{(i,j)\ ,(i,k)\ ,(k,i)\ ,(k,j)\ ,(j,k)\}=E_H
$$
cannot occur because $E_+|_{{\overline{H}}}=\{\{i,k\},\{k,j\}\}$ and $E_-|_{{\overline{H}}}=\emptyset$ 
are not permitted by the signed eliminability. Taking 
into account the fact that $\overline{G}$ is 
signed eliminable, using symmetry, we have two possibilities of $H$:

\begin{enumerate}
\item
$\{(i,j),\ (j,i),\ (i,k),\ (k,i),\ (j,k) \}= E_H$.
\item
$\{(i,j),\ (j,i),\ (i,k),\ (k,i),\ (k,j) \}= E_H$.
\end{enumerate}

Then both cases satisfy (A1) and (A2).
\medskip

\textbf{Case 5}. $E_+|_{{\overline{H}}}=\emptyset$ and $|E_-|_{{\overline{H}}}|=2$. 
In this case, taking into account the fact that $\overline{G}$ is 
signed eliminable, using symmetry, we have two possibilities of $H$:

\begin{enumerate}
\item
$\{(i,k) \}= E_H$.
\item
$\{(k,i) \}= E_H$.
\end{enumerate}

Then both cases satisfy (A1) and (A2).
\medskip

\textbf{Case 6}. $|E_+|_{{\overline{H}}}|=1$ and $|E_-|_{{\overline{H}}}|=1$. 
In this case, taking into account the fact that $\overline{G}$ is 
signed eliminable, using symmetry, we have two possibilities of $H$:

\begin{enumerate}
\item
$\{(i,j),\ (i,k),\ (k,i) \}=E_H$.
\item
$\{(j,i),\ (i,k),\ (k,i) \}=E_H$.
\end{enumerate}

Then both cases satisfy (A1) and (A2).
\medskip

\textbf{Case 7}. $|E_+|_{{\overline{H}}}|=3$ and $E_-|_{{\overline{H}}}=\emptyset$. 
In this case, we have only one possibilities of $H$:

\begin{enumerate}
\item
$\{(i,j),\ (j,i),\ (k,i),\ (i,k),\ (j,k),\ (k,j) \}= E_H$.
\end{enumerate}

This digraph satisfies (A1) and (A2).
\medskip

\textbf{Case 8}. $E_+|_{{\overline{H}}}=\emptyset$ and $|E_-|_{{\overline{H}}}|=3$. 
In this case, $H$ has no edges. Hence it satisfies (A1) and (A2). 
\medskip

\textbf{Case 9}. $|E_+|_{{\overline{H}}}|=2$ and $|E_-|_{{\overline{H}}}|=1$. 
In this case, taking into account the fact that $\overline{G}$ is 
signed eliminable, using symmetry, we have only one possibilities of $H$:

\begin{enumerate}
\item
$\{(i,j),\ (j,i),\ (k,i),\ (i,k) \}= E_H$.
\end{enumerate}

This satisfies (A1) and (A2).
\medskip

\textbf{Case 10}. $|E_+|_{{\overline{H}}}|=1$ and $|E_-|_{{\overline{H}}}|=2$. 
In this case, taking into account the fact that $\overline{G}$ is 
signed eliminable, using symmetry, we have only one possibilities of $H$:

\begin{enumerate}
\item
$\{(i,k),\ (k,i) \}= E_H$.
\end{enumerate}

This satisfies (A1) and (A2).
\medskip

The classification above shows that the signed elimination ordering of 
$\overline{G}$ makes $G$ into a digraph satisfying (A1) and (A2), which completes the proof. \owari
\medskip

%The next corollary is a characterization of a digraph satisfying (A1) and (A2) without using 
%the numbering of vertices. The proof follows immediately by the above arguments.

%If $H$ satisfies the conditions (A1) and (A2), then 
%there is noting to show. The only cases when $H$ doses not satisfy the two conditions
%are the following:
%
%
%
%However, these three $H$'s cannot admit any numbering of 
%vertices which makes $H$ into A-graphs. Also, it is an easy computation to show that 
%these three $H$'s produce deformations of braid arrangement the conings of which are 
%not free, see \cite{A1} for example. In other words, it is easy to check a conjecture of 
%Athanasiadis when $\ell=2$. So if (1), (2) or (3) happens, then 
%the localization of $c\A_G$ at $X=\{x_i=x_j=x_k,z\} \in L(c\A_G)$ 
%is not free, which contradicts the freeness of $c\A_G$. Hence 
%these three cases never happen, which proves nothing but the signed eliminable ordering 
%of $V_{\overline{G}}$ induces a numbering of $V_G$ which makes $G$ into an A-graph. 
%Hence if $c\A_G$ is free, then $G$ is an A-graph, which completes the proof. \owari

\noindent
\textbf{Proof of Theorem \ref{nece}}.
Let $G=(V_G,E_G)$ be a digraph and 
$\A_G$ the corresponding deformation of a braid arrangement of the form (\ref{defath}). 
Assume that $c\A_G$ is free. Take the Ziegler restriction of $c\A_G$ 
onto the infinite hyperplane $H_\infty \in c\A_G$. Then we obtain a multiarrangement 
$(\A_\ell,m_{\overline{G}})$, where $\A_\ell$ is a braid arrangement in $\K^{\ell+1}$, 
$S(G)=\overline{G}=(V_G,E_+ \cup E_-)$ is a signed graph defined above and 
$m_{\overline{G}}$ is the corresponding multiplicity to $\overline{G}$. See \cite{ANN} for details.
%Note that $G$ and $\overline{G}$ have the same set of 
%vertices. 
%By this correspondence we can construct a multiplicity 
%$m_{\overline{G}}:\A_\ell \rightarrow \Z_{>0}$ and a multiarrangement 
%$(\A_\ell,m_{\overline{G}})$. See 
%\cite{ANN} for details. 
By \cite{Z}, the Ziegler restriction of a free central arrangement 
is a free multiarrangement. Hence $(\A_\ell,m_{\overline{G}})$ is free. 
Then Theorem 0.3 in \cite{ANN} shows that $\overline{G}$ is signed eliminable. 
For any codimension three flat $X \in L(c\A_G)$ with $X \subset H_\infty$, 
the localization $(c\A_G)_X$ is free since $c\A_G$ is free. Let 
$X=\{x_i=x_j=x_k\} \cap H_\infty$. Then 
$(c\A_G)_X$ is a deformation of $A_2$-type arrangements. 
Also, by definition of the deformation (\ref{defath}), there is an 
induced subgraph $H$ of $G$ with $V_H=\{i,j,k\}$ such that 
$(c\A_G)_X=c\A_H \times \emptyset_{\ell-2}$.
%, where $\emptyset_{\ell-2}$ is an 
%empty arrangement in $\K^{\ell-2}$. 
%Recalling the fact that the freeness of 
%$c\A_H \times \emptyset_{\ell-2}$ is equivalent to that of $c\A_H$, 
%Hence the deformation $(c\A_G)_X=c\A_H$ corresponds to some induced subgraph of $G$ with 
%three vertices. 
Then freeness of $(c\A_G)_X=c\A_H \times \emptyset_{\ell-2}$ and Lemma \ref{lemma1} show that 
the subgraphs (1), (2) and (3) in Lemma \ref{lemma1} are not contained in $G$. Then Proposition \ref{char} 
implies that $G$ satisfies (A1) and (A2). \owari

% \vspace{5mm}

%\noindent
%Takuro Abe\\
%Department of Mathematics\\
%Hokkaido University\\
%Sapporo 060--0810, Japan\\
%abetaku@math.sci.hokudai.ac.jp
%
%\bigskip
%
%\noindent
%Koji Nuida\\
%Research Center for Information Security (RCIS)\\
%National Institute of Advanced Industrial Science and Technology (AIST)\\
%Tokyo 101--0021, Japan\\
%k.nuida@aist.go.jp
%
%%\noindent%
%
%%{\em e-mail address}\ : \ 
%
%\bigskip
%
%\noindent
%Yasuhide Numata\\
%Department of Mathematics\\
%Hokkaido University\\
%Sapporo 060--0810, Japan\\
%nu@math.sci.hokudai.ac.jp

\end{document}